%% file: 2021lindenhovius-j-to-X.tex
\documentclass[oneside,11pt]{article}
\usepackage[colorlinks,citecolor=DarkRed,urlcolor=DarkRed]{hyperref} % (find-es "tex" "hyperref")
\usepackage{amsmath}
\usepackage{amsfonts}
\usepackage{amssymb}
\usepackage{pict2e}
\usepackage[x11names,svgnames]{xcolor} % (find-es "tex" "xcolor")
\usepackage{ifluatex}
\ifluatex
  \input edrxaccents.tex            % (find-LATEX "edrxaccents.tex")
  \input edrx21chars.tex            % (find-LATEX "edrxchars.tex")
  \input edrxheadfoot.tex           % (find-LATEX "edrxheadfoot.tex")
\else
  \usepackage[utf8]{inputenc}
  \input edrx21chars-d.tex          % (find-LATEX "edrx21chars-d.tex")
\fi
\usepackage{edrx21}                 % (find-LATEX "edrx21.sty")
\usepackage{proof}   % For derivation trees ("%:" lines)
\input diagxy        % For 2D diagrams ("%D" lines)
\usepackage[backend=biber,
   style=alphabetic]{biblatex}            % (find-es "tex" "biber")
\begin{document}

\ifluatex
  \catcode`\^^J=10
  \directlua{dofile "dednat6load.lua"}
\else
  \input\jobname.dnt   % (find-LATEXfile "2021lindenhovius-j-to-X.dnt")
  \def\pu{}
\fi

%L forths["<.>"]  = function () pusharrow("<.>") end
%L forths["<-->"] = function () pusharrow("<-->") end
%L forths["|-->"] = function () pusharrow("|-->") end
%L forths["<--|"] = function () pusharrow("<--|") end

% «defs»  (to ".defs")

\def\Downs {\mathcal{D}}
\def\Ddp   {\Downs({↓}p)}
\def\singp {\{p\}}
\def\dnp   {{↓}p}
\def\dnpo  {{↓}p∖\{p\}}
\def\setofsc#1#2{\{\,#1\;:\;#2\,\}}
\def\nuc   {(·)^*}
\def\Nucs  {\mathsf{Nucs}}
\def\GrTops{\mathsf{GrTops}}

%  _____ _ _   _      
% |_   _(_) |_| | ___ 
%   | | | | __| |/ _ \
%   | | | | |_| |  __/
%   |_| |_|\__|_|\___|
%                     
% «title»  (to ".title")
% (linfp 1 "title")
% (linfa   "title")

\title{On a formula that is not in \\ ``Grothendieck Topologies in
  Posets''}

\author{
  % A.J.~Lindenhovius
  % \and
  Eduardo Ochs
  }

\maketitle

\begin{abstract}

  The paper \cite{Lindenhovius} shows that when $𝐏$ is an Artinian
  poset and $𝐄$ is the topos $\Set^𝐏$ then there are bijections
  between the set of subsets of $𝐏$, the set of Grothendieck
  topologies on $𝐄$, and the set of nuclei on the Heyting Algebra
  $\Sub(1_𝐄)$. It also shows that there are nice formulas for
  converting between subsets, Grothendieck topologies, and nuclei, but
  the formula for converting a nucleus to a subset is not spelled out
  explicitly. These notes fix that gap.

  % The paper ``Grothendieck Topologies on Posets'' by A.J.
  % Lindenhovius shows that when $\mathbf{P}$ is an Artinian poset and
  % $\mathbf{E}$ is the topos $\mathbf{Set}^\mathbf{P}$ then there are
  % bijections between the set of subsets of $\mathbf{P}$, the set of
  % Grothendieck topologies on $\mathbf{E}$, and the set of nuclei on
  % the Heyting Algebra $\mathrm{Sub}(1_\mathbf{E})$. It also shows
  % that there are nice formulas for converting between subsets,
  % Grothendieck topologies, and nuclei, but the formula for
  % converting a nucleus to a subset is not spelled out explicitly.
  % These notes fix that gap.

\end{abstract}

Let $𝐏$ be a downward-directed poset, that we will also regard as a
category. Then $𝐄:=\Set^𝐏$ is a topos and $H:=\Sub(1_𝐄)$ is a Heyting
Algebra. Let's denote the set of all Grothendieck topologies on $𝐄$ by
$\GrTops(𝐏)$, the set of all nuclei on $H$ by $\Nucs(𝐏)$, and the set
of subsets of (the set of points of) $𝐏$ by $\Pts(𝐏)$. In
\cite{Lindenhovius} Bert Lindenhovius shows that when $𝐏$ is Artinian
we have bijections between $\Pts(𝐏)$, $\Nucs(𝐏)$, and $\GrTops(𝐏)$. He
uses notational conventions in which $J$ always denotes a Grothendieck
topology, $j$ always denotes a nucleus, and $X$ always denotes a
subset of $𝐏$, and he writes the components these bijections as
$(j↦J_j)$, $(J↦j_J)$, and so on; we will also write them here as
$(j↦J)$, $(J↦j)$, etc, and we will write the bijections as $(X↔j)$,
$(X↔J)$, and $(j↔J)$. Let's put all this in a diagram:
%
%D diagram bijections
%D 2Dx     100  +45 +40  +45
%D 2D  100 A0 - A1  B0 - B1
%D 2D      |  /     |  /
%D 2D  +45 A2 - A3  B2 - B3
%D 2D
%D ren A0 A1 A2 ==> \Pts(𝐏) \Nucs(𝐏) \GrTops(𝐏)
%D ren B0 B1 B2 ==> X j J
%D
%D (( A0 A1  -> sl^ .plabel= a \smt{C.4.2}
%D    A0 A1 <-- sl_ .plabel= b \smt{(e-mail)}
%D    A0 A2  -> sl_ .plabel= l \smtt{2.8,}{C.4.1}
%D    A0 A2 <- sl^ .plabel= r \smt{2.9}
%D    A2 A1 <-> .plabel= m \smtt{B.8,}{B.25}
%D    # A2 A3 <-->
%D ))
%D (( B0 B1  |-> sl^ .plabel= a \sm{(X↦j)}
%D    B0 B1 <--| sl_ .plabel= b \sm{(j↦X)}
%D    B0 B2  <->     .plabel= l \sm{(X↦J),\\(J↦X)}
%D    B2 B1  <->     .plabel= m \sm{(j↦J),\\(J↦j)}
%D ))
%D enddiagram
%D
$$\pu
  \diag{bijections}
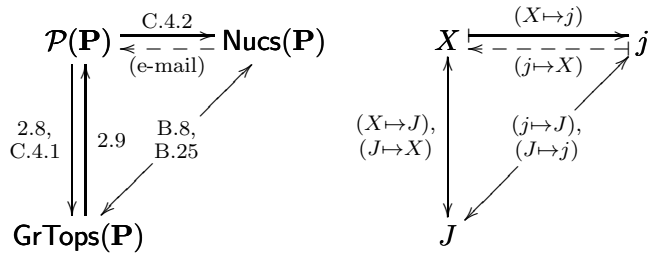
$$

\newpage

% «components»  (to ".components")
% (linfp 2 "components")
% (linfa   "components")

He defines the components of these bijections as:
$$\begin{array}{rrcll}
  (X↦j): & j_X(S) &=&  X→S                               & \text{(C.4.2, C.2)} \\
  (j↦X): & X_j    &=& \setofsc {p∈𝐏} {p\not∈j(\dnpo)}    & \text{(e-mail)} \\
  [5pt]
  (X↦J): & J_X(p) &=& \setofsc{S∈\Ddp}{X∩{↓}p⊆S}         & \text{(2.8, C.4.1)} \\
  (J↦X): & X_J    &=& \setofsc{p∈𝐏}{J(p)=\{{↓}p\}}      & \text{(2.9)} \\
  [5pt]
  (j↦J): & J_j(p) &=& \setofsc {S∈\Downs(\dnp)} {p∈j(S)} & \text{(B.8, B.25)} \\
  (J↦j): & j_J(S) &=& \setofsc {p∈𝐏} {S∩{↓}p∈J(p)}      & \text{(B.8, B.25)} \\
  \end{array}
$$

The annotations like ``(C.4.2, C.2)'' indicate where these components
are defined. Note that one of the annotations says ``(e-mail)''; this
is because that formula doesn't appear explicitly in
\cite{Lindenhovius}, and so I (Eduardo) asked him (Bert) if that
formula was what I guessed it would be, and he replied with a formula
slightly shorter than my guess, and a proof...

\msk

These notes are just to make his formula and his proof available in a
public place. All the mathematical content here is by Bert
Lindenhovius, and all the typesetting was done by Eduardo Ochs, who
found Bert's proof hard to follow and decided to typeset it in Natural
Deduction form. {\sl Note:} when I wrote the first version of these
notes I listed Bert as the main author and me as the coauthor, but he
told me that he preferred to be credited only in the text. I insisted,
and explained that I would have never been able to find these proofs
by myself --- but he insisted more.

\bsk

If we combine $(j↦J)$ and $(J↦X)$ we get this:
$$\begin{array}{crcl}
    (j↦J): & J_j(p)  &=& \setofsc {S∈\Downs(\dnp)} {p∈j(S)} \\
    (J↦X): & X_J     &=& \setofsc {p∈𝐏} {J(p)=\{{↓}p\}}       \\
  (j↦J↦X): & X_{J_j} &=& \setofsc {p∈𝐏} { \setofsc {S∈\Downs(\dnp)} {p∈j(S)} =\{{↓}p\}}       \\
           &         &=& \setofsc {p∈𝐏} {∀S∈\Downs(\dnp). \; ((p∈j(S))↔(S={↓}p))}       \\
        [5pt]
    (j↦X): & X_j     &=& \setofsc {p∈𝐏} {p\not∈j(\dnpo)} \\
  \end{array}
$$

It is not obvious at all that $X_{J_j} = X_j$. We will prove that
$p∈X_j$ iff $p∈X_{J_j}$, where:
$$\begin{array}{rcl}
    (p∈X_j)     &=& (p\not∈j(\dnpo)) \\
    (p∈X_{J_j}) &=& ∀S∈\Downs(\dnp). \; (p∈j(S))↔(S={↓}p) \\
  \end{array}
$$

\newpage

% «proofs»  (to ".proofs")
% (linfp 3 "proofs")
% (linfa   "proofs")

\def\H#1{\hspace{#1cm}}

Look:
%:
%:
%:                                        p∈X_{J_j}
%:   ------------------   ---------------------------------
%:   \dnpo∈\Downs(\dnp)   ∀S∈\Downs(\dnp).(p∈j(S))↔(S=\dnp)
%:   ------------------------------------------------------
%:       (p∈j(\dnpo))↔(\dnpo=\dnp)
%:       ---------------------------------     --------------
%:       (p\not∈j(\dnpo))↔(\dnpo\not=\dnp)     \dnpo\not=\dnp
%:       ----------------------------------------------------
%:       p\not∈j(\dnpo)
%:       --------------------
%:       p∈X_j
%:
%:       ^part-small
%:
\pu
$$\ded{part-small}$$
%:
%:
%:   [S∈\Downs(\dnp)]^2
%:   ------------------
%:     S⊂\dnp            [S\neq\dnp]^1     [S∈\Downs(\dnp)]^2
%:     -------------------------------     --------------
%:          p\not∈S                        S⊂\dnp
%:          -------------------------------------
%:                  S⊂\dnpo                                p∈X_j
%:               -------------                         --------------
%:               j(S)⊂j(\dnpo)             \H{-1}      p\not∈j(\dnpo)
%:               ----------------------------------------------------
%:                  p\not∈j(S)
%:                  ------------------------1
%:                  (S\neq\dnp)→(p\not∈j(S))
%:                  ------------------------            -----------------
%:                  (p∈j(S))→(S=\dnp)          \H{-1}   (S=\dnp)→(p∈j(S))
%:                  -----------------------------------------------------
%:                  (p∈j(S))↔(S=\dnp)
%:                  ---------------------------------2
%:                  ∀S∈\Downs(\dnp).(p∈j(S))↔(S=\dnp)
%:                  ---------------------------------
%:                  p∈X_{J_j}
%:
%:                  ^part-big
%:
\pu
$$\ded{part-big}$$

This proves that $X_j$ and $X_{J_j}$ are equal.

\newpage

% «other-formula»  (to ".other-formula")
% (linfp 4 "other-formula")
% (linfa   "other-formula")

Now let's check that $X_j$ and $X'_j$ are equal,

where $X'_j$ is defined as:

$$\begin{array}{rcl}
    X_j  &=& \setofsc {p∈𝐏} {p\not∈j(\dnpo)} \\
    X'_j &=& \setofsc {p∈𝐏} {j(\dnp)\not=j(\dnpo)} \\
  \end{array}
$$

Proof:
%:
%:
%:   
%:   [p∈j(\dnpo)]^1
%:   --------------            ----------         ------------
%:   \dnp⊆j(\dnpo)             \dnpo⊆\dnp         \dnp⊆j(\dnp)   [j(\dnp)=j(\dnpo)]^1
%:   --------------------    ----------------     -----------------------------------
%:   j(\dnp)⊆(j∘j)(\dnpo)    j(\dnpo)⊆j(\dnp)     \dnp⊆j(\dnp)=j(\dnpo)       
%:   ----------------------------------------     ---------------------       
%:   j(\dnp)⊆(j∘j)(\dnpo)=j(\dnpo)⊆j(\dnp)        \dnp⊆j(\dnpo)               
%:   -------------------------------------        -------------               
%:   j(\dnp)=j(\dnpo)                                p∈j(\dnpo)                  
%:   ----------------------------------------------------------1
%:   (p∈j(\dnpo))↔(j(\dnp)=j(\dnpo))
%:   ---------------------------------------
%:   (p\not∈j(\dnpo))↔(j(\dnp)\not=j(\dnpo))
%:   ---------------------------------------
%:   p∈X_j↔p∈X'_j
%:
%:   ^lemma
%:
%: 
\pu
%
% (lindp 50 "B.6")
% (linda    "B.6")
%
$$\ded{lemma}$$

The formula of $X'_j$ above is the one that I asked if it was correct;
Bert answered that yes, and showed that it is equivalent to the
slighty shorter formula for $X_j$.

%L write_dnt_file()
\pu

\printbibliography

\GenericWarning{Success:}{Success!!!}  % Used by `M-x cv'

\end{document}

%% file: edrx21chars-d.tex
\DeclareUnicodeCharacter{00A1}{\text{\textexclamdown}} % ¡
\DeclareUnicodeCharacter{00A7}{\S}                 % §
\DeclareUnicodeCharacter{00AC}{\neg}               % ¬
\DeclareUnicodeCharacter{00B0}{^\circ}             % °
\DeclareUnicodeCharacter{00B1}{\pm}                % ±
\DeclareUnicodeCharacter{00B2}{^2}                 % ²
\DeclareUnicodeCharacter{00B3}{^3}                 % ³
\DeclareUnicodeCharacter{00B7}{\cdot}              % ·
\DeclareUnicodeCharacter{00B9}{^{-1}}              % ¹
\DeclareUnicodeCharacter{00D7}{\times}             % ×
\DeclareUnicodeCharacter{00F7}{\div}               % ÷
\DeclareUnicodeCharacter{0393}{\Gamma}             % Γ
\DeclareUnicodeCharacter{0394}{\Delta}             % Δ
\DeclareUnicodeCharacter{0398}{\Theta}             % Θ
\DeclareUnicodeCharacter{039B}{\Lambda}            % Λ
\DeclareUnicodeCharacter{03A0}{\Pi}                % Π
\DeclareUnicodeCharacter{03A3}{\Sigma}             % Σ
\DeclareUnicodeCharacter{03A6}{\Phi}               % Φ
\DeclareUnicodeCharacter{03A8}{\Psi}               % Ψ
\DeclareUnicodeCharacter{03A9}{\Omega}             % Ω
\DeclareUnicodeCharacter{03B1}{\alpha}             % α
\DeclareUnicodeCharacter{03B2}{\beta}              % β
\DeclareUnicodeCharacter{03B3}{\gamma}             % γ
\DeclareUnicodeCharacter{03B4}{\delta}             % δ
\DeclareUnicodeCharacter{03B5}{\epsilon}           % ε
\DeclareUnicodeCharacter{03B6}{\zeta}              % ζ
\DeclareUnicodeCharacter{03B7}{\eta}               % η
\DeclareUnicodeCharacter{03B8}{\theta}             % θ
\DeclareUnicodeCharacter{03B9}{\iota}              % ι
\DeclareUnicodeCharacter{03BA}{\kappa}             % κ
\DeclareUnicodeCharacter{03BB}{\lambda}            % λ
\DeclareUnicodeCharacter{03BC}{\mu}                % μ
\DeclareUnicodeCharacter{03BD}{\nu}                % ν
\DeclareUnicodeCharacter{03BE}{\xi}                % ξ
\DeclareUnicodeCharacter{03C0}{\pi}                % π
\DeclareUnicodeCharacter{03C1}{\rho}               % ρ
\DeclareUnicodeCharacter{03C3}{\sigma}             % σ
\DeclareUnicodeCharacter{03C4}{\tau}               % τ
\DeclareUnicodeCharacter{03C6}{\phi}               % φ
\DeclareUnicodeCharacter{03C7}{\chi}               % χ
\DeclareUnicodeCharacter{03C8}{\psi}               % ψ
\DeclareUnicodeCharacter{03C9}{\omega}             % ω
\DeclareUnicodeCharacter{03D5}{\origphi}           % ϕ
\DeclareUnicodeCharacter{1D400}{\catA}             % 𝐀
\DeclareUnicodeCharacter{1D401}{\catB}             % 𝐁
\DeclareUnicodeCharacter{1D402}{\catC}             % 𝐂
\DeclareUnicodeCharacter{1D403}{\catD}             % 𝐃
\DeclareUnicodeCharacter{1D404}{\catE}             % 𝐄
\DeclareUnicodeCharacter{1D405}{\catF}             % 𝐅
\DeclareUnicodeCharacter{1D406}{\catG}             % 𝐆
\DeclareUnicodeCharacter{1D407}{\catH}             % 𝐇
\DeclareUnicodeCharacter{1D408}{\catI}             % 𝐈
\DeclareUnicodeCharacter{1D409}{\catJ}             % 𝐉
\DeclareUnicodeCharacter{1D40A}{\catK}             % 𝐊
\DeclareUnicodeCharacter{1D40B}{\catL}             % 𝐋
\DeclareUnicodeCharacter{1D40C}{\catM}             % 𝐌
\DeclareUnicodeCharacter{1D40D}{\catN}             % 𝐍
\DeclareUnicodeCharacter{1D40E}{\catO}             % 𝐎
\DeclareUnicodeCharacter{1D40F}{\catP}             % 𝐏
\DeclareUnicodeCharacter{1D410}{\catQ}             % 𝐐
\DeclareUnicodeCharacter{1D411}{\catR}             % 𝐑
\DeclareUnicodeCharacter{1D412}{\catS}             % 𝐒
\DeclareUnicodeCharacter{1D413}{\catT}             % 𝐓
\DeclareUnicodeCharacter{1D414}{\catU}             % 𝐔
\DeclareUnicodeCharacter{1D415}{\catV}             % 𝐕
\DeclareUnicodeCharacter{1D416}{\catW}             % 𝐖
\DeclareUnicodeCharacter{1D417}{\catX}             % 𝐗
\DeclareUnicodeCharacter{1D418}{\catY}             % 𝐘
\DeclareUnicodeCharacter{1D419}{\catZ}             % 𝐙
\DeclareUnicodeCharacter{1D41B}{\mathbf}           % 𝐛
\DeclareUnicodeCharacter{1D422}{\textsl}           % 𝐢
\DeclareUnicodeCharacter{1D42B}{\mathrm}           % 𝐫
\DeclareUnicodeCharacter{1D42C}{\mathsf}           % 𝐬
\DeclareUnicodeCharacter{1D42D}{\text}             % 𝐭
\DeclareUnicodeCharacter{1D4D0}{\mathcal{A}}       % 𝓐
\DeclareUnicodeCharacter{1D4D1}{\mathcal{B}}       % 𝓑
\DeclareUnicodeCharacter{1D4D2}{\mathcal{C}}       % 𝓒
\DeclareUnicodeCharacter{1D4D3}{\mathcal{D}}       % 𝓓
\DeclareUnicodeCharacter{1D4D4}{\mathcal{E}}       % 𝓔
\DeclareUnicodeCharacter{1D4D5}{\mathcal{F}}       % 𝓕
\DeclareUnicodeCharacter{1D4D6}{\mathcal{G}}       % 𝓖
\DeclareUnicodeCharacter{1D4D7}{\mathcal{H}}       % 𝓗
\DeclareUnicodeCharacter{1D4D8}{\mathcal{I}}       % 𝓘
\DeclareUnicodeCharacter{1D4D9}{\mathcal{J}}       % 𝓙
\DeclareUnicodeCharacter{1D4DA}{\mathcal{K}}       % 𝓚
\DeclareUnicodeCharacter{1D4DB}{\mathcal{L}}       % 𝓛
\DeclareUnicodeCharacter{1D4DC}{\mathcal{M}}       % 𝓜
\DeclareUnicodeCharacter{1D4DD}{\mathcal{N}}       % 𝓝
\DeclareUnicodeCharacter{1D4DE}{\mathcal{O}}       % 𝓞
\DeclareUnicodeCharacter{1D4DF}{\mathcal{P}}       % 𝓟
\DeclareUnicodeCharacter{1D4E0}{\mathcal{Q}}       % 𝓠
\DeclareUnicodeCharacter{1D4E1}{\mathcal{R}}       % 𝓡
\DeclareUnicodeCharacter{1D4E2}{\mathcal{S}}       % 𝓢
\DeclareUnicodeCharacter{1D4E3}{\mathcal{T}}       % 𝓣
\DeclareUnicodeCharacter{1D4E4}{\mathcal{U}}       % 𝓤
\DeclareUnicodeCharacter{1D4E5}{\mathcal{V}}       % 𝓥
\DeclareUnicodeCharacter{1D4E6}{\mathcal{W}}       % 𝓦
\DeclareUnicodeCharacter{1D4E7}{\mathcal{X}}       % 𝓧
\DeclareUnicodeCharacter{1D4E8}{\mathcal{Y}}       % 𝓨
\DeclareUnicodeCharacter{1D4E9}{\mathcal{Z}}       % 𝓩
\DeclareUnicodeCharacter{2022}{\bullet}            % •
\DeclareUnicodeCharacter{2026}{\ldots}             % …
\DeclareUnicodeCharacter{2113}{\ell}               % ℓ
\DeclareUnicodeCharacter{214B}{\bindnasrepma}      % ⅋
\DeclareUnicodeCharacter{2190}{\ot}                % ←
\DeclareUnicodeCharacter{2191}{\upto}              % ↑
\DeclareUnicodeCharacter{2192}{\to}                % →
\DeclareUnicodeCharacter{2193}{\dnto}              % ↓
\DeclareUnicodeCharacter{2194}{\bij}               % ↔
\DeclareUnicodeCharacter{2195}{\updownarrow}       % ↕
\DeclareUnicodeCharacter{2196}{\nwarrow}           % ↖
\DeclareUnicodeCharacter{2197}{\nearrow}           % ↗
\DeclareUnicodeCharacter{2198}{\searrow}           % ↘
\DeclareUnicodeCharacter{2199}{\swarrow}           % ↙
\DeclareUnicodeCharacter{21A3}{\twoheadrightarrow} % ↣
\DeclareUnicodeCharacter{21A4}{\mapsfrom}          % ↤
\DeclareUnicodeCharacter{21A6}{\mapsto}            % ↦
\DeclareUnicodeCharacter{21C0}{\rightharpoonup}    % ⇀
\DeclareUnicodeCharacter{21D0}{\Leftarrow}         % ⇐
\DeclareUnicodeCharacter{21D2}{\funto}             % ⇒
\DeclareUnicodeCharacter{21D4}{\Leftrightarrow}    % ⇔
\DeclareUnicodeCharacter{2200}{\forall}            % ∀
\DeclareUnicodeCharacter{2202}{\partial}           % ∂
\DeclareUnicodeCharacter{2203}{\exists}            % ∃
\DeclareUnicodeCharacter{2205}{\emptyset}          % ∅
\DeclareUnicodeCharacter{2207}{\nabla}             % ∇
\DeclareUnicodeCharacter{2208}{\in}                % ∈
\DeclareUnicodeCharacter{2216}{\backslash}         % ∖
\DeclareUnicodeCharacter{2218}{\circ}              % ∘
\DeclareUnicodeCharacter{221A}{\sqrt}              % √
\DeclareUnicodeCharacter{221E}{\infty}             % ∞
\DeclareUnicodeCharacter{2227}{\land}              % ∧
\DeclareUnicodeCharacter{2228}{\lor}               % ∨
\DeclareUnicodeCharacter{2229}{\cap}               % ∩
\DeclareUnicodeCharacter{222A}{\cup}               % ∪
\DeclareUnicodeCharacter{222B}{\int}               % ∫
\DeclareUnicodeCharacter{223C}{\sim}               % ∼
\DeclareUnicodeCharacter{2243}{\simeq}             % ≃
\DeclareUnicodeCharacter{2245}{\cong}              % ≅
\DeclareUnicodeCharacter{2248}{\approx}            % ≈
\DeclareUnicodeCharacter{2260}{\neq}               % ≠
\DeclareUnicodeCharacter{2261}{\equiv}             % ≡
\DeclareUnicodeCharacter{2264}{\le}                % ≤
\DeclareUnicodeCharacter{2265}{\ge}                % ≥
\DeclareUnicodeCharacter{2282}{\subset}            % ⊂
\DeclareUnicodeCharacter{2283}{\supset}            % ⊃
\DeclareUnicodeCharacter{2286}{\subseteq}          % ⊆
\DeclareUnicodeCharacter{2287}{\supseteq}          % ⊇
\DeclareUnicodeCharacter{2293}{\sqcap}             % ⊓
\DeclareUnicodeCharacter{2294}{\sqcup}             % ⊔
\DeclareUnicodeCharacter{2295}{\oplus}             % ⊕
\DeclareUnicodeCharacter{2296}{\ominus}            % ⊖
\DeclareUnicodeCharacter{2297}{\otimes}            % ⊗
\DeclareUnicodeCharacter{2298}{\oslash}            % ⊘
\DeclareUnicodeCharacter{2299}{\odot}              % ⊙
\DeclareUnicodeCharacter{22A2}{\vdash}             % ⊢
\DeclareUnicodeCharacter{22A3}{\dashv}             % ⊣
\DeclareUnicodeCharacter{22A4}{\top}               % ⊤
\DeclareUnicodeCharacter{22A5}{\bot}               % ⊥
\DeclareUnicodeCharacter{22A8}{\vDash}             % ⊨
\DeclareUnicodeCharacter{22C0}{\bigwedge}          % ⋀
\DeclareUnicodeCharacter{22C1}{\bigvee}            % ⋁
\DeclareUnicodeCharacter{22C2}{\bigcap}            % ⋂
\DeclareUnicodeCharacter{22C3}{\bigcup}            % ⋃
\DeclareUnicodeCharacter{22C4}{\lozenge}           % ⋄
\DeclareUnicodeCharacter{22C5}{\Box}               % ⋅
\DeclareUnicodeCharacter{2329}{\langle}            % 〈
\DeclareUnicodeCharacter{232A}{\rangle}            % 〉
\DeclareUnicodeCharacter{2581}{\_}                 % ▁
\DeclareUnicodeCharacter{25A1}{\Box}               % □
\DeclareUnicodeCharacter{25FB}{\Box}               % ◻
\DeclareUnicodeCharacter{266D}{\flat}              % ♭
\DeclareUnicodeCharacter{266E}{\natural}           % ♮
\DeclareUnicodeCharacter{266F}{\sharp}             % ♯
\DeclareUnicodeCharacter{27E6}{\llbracket}         % ⟦
\DeclareUnicodeCharacter{27E7}{\rrbracket}         % ⟧
\DeclareUnicodeCharacter{2806}{{:}}                % ⠆